\documentclass[12pt,leqno]{article}

\usepackage[active]{srcltx}

\usepackage{amsfonts,latexsym,amsmath,amssymb,amsthm}
\usepackage{fullpage}
\newtheorem{theorem}{Theorem}[section]
\newtheorem{lemma}[theorem]{Lemma}

\newtheorem{corollary}[theorem]{Corollary}
\newtheorem{remark}[theorem]{Remark}
\newtheorem{question}[theorem]{Question}

\theoremstyle{definition}

\def\R{{\mathbb R}}
\def\e{\varepsilon}
\def\N{{\mathbb N}}

\theoremstyle{remark}

\newtheorem*{note*}{Note}
\numberwithin{equation}{section}
\makeatletter
\newcommand{\rank}{\mathop{\operator@font rank}}
\newcommand{\conv}{\mathop{\operator@font conv}}
\newcommand{\vol}{\mathop{\operator@font vol}}
\newcommand{\onetagright}{\tagsleft@false}
\makeatother
\newcommand{\ls}{\leqslant}
\newcommand{\gr}{\geqslant}
\renewcommand{\epsilon}{\varepsilon}

\usepackage{color}

\usepackage{hyperref}

\begin{document}
\small

\title{\bf Volume difference inequalities}

\medskip

\author{Apostolos Giannopoulos and Alexander Koldobsky}

\date{}

\maketitle
\begin{abstract} We prove several inequalities estimating the distance between volumes of two
bodies in terms of the maximal or minimal difference between areas of sections or projections
of these bodies. We also provide extensions in which volume is replaced by an arbitrary measure.
\end{abstract}

%%%%%%%%%%%%%%%%%%%%%%%%%%%%%%%%%%%%%%%%%%%%%%%%%%%%%%%%%%%%%%%%%%%%%%%%%%%%%%%%%%%%%%%%%%%%%%%%%%%%%%%%%%%%%%%%%%%%%%%%%%%%%%%%%%%%%%%%%%%%%%%%%%
\section{Introduction}
%%%%%%%%%%%%%%%%%%%%%%%%%%%%%%%%%%%%%%%%%%%%%%%%%%%%%%%%%%%%%%%%%%%%%%%%%%%%%%%%%%%%%%%%%%%%%%%%%%%%%%%%%%%%%%%%%%%%%%%%%%%%%%%%%%%%%%%%%%%%%%%%%%
Volume difference inequalities are designed to estimate the error in computations of volume of a body out of the areas of its sections and projections. We start with the case of sections. Let $\gamma_{n,k}$ be the smallest constant $\gamma>0$ satisfying the inequality
\begin{equation} \label{main-prob1}
|K|^{\frac {n-k}n}-|L|^{\frac {n-k}n}\ls \gamma^k
\max_{F\in {\rm Gr}_{n-k}} \big (|K\cap F|-|L\cap F|\big )
\end{equation}
for all $1\ls k <n$ and all origin-symmetric convex bodies $K$ and $L$ in $\R^n$ such that $L\subset K.$
Here ${\rm Gr}_{n-k}$ is the Grassmanian of $(n-k)$-dimensional subspaces of $\R^n,$ and
$|K|$ stands for volume of appropriate dimension.

\begin{question} \label{q-sect}
Does there exist an absolute constant $C$ so that $\sup_{n,k} \gamma_{n,k}\ls C\ ?$
\end{question}

Question \ref{q-sect} is stronger than the slicing problem, a major open problem in convex geometry \cite{Bourgain-1986, Bourgain-1987, Ball-thesis, MP}. In fact,
putting $L=\beta B_2^n$ in \eqref{main-prob1},
where $B_2^n$ is the unit Euclidean ball in $\R^n,$ and then sending $\beta$ to zero,
one gets the slicing problem: does there exist an absolute
constant $C$ so that for any $1\ls k <n,$ and any origin-symmetric convex body $K$ in $\R^n$
\begin{equation} \label{slicing}
|K|^{\frac {n-k}n}\ls C^k
\max_{H\in {\rm Gr}_{n-k}} |K\cap H|\ ?
\end{equation}
The best-to-date general estimate $C\ls O(n^{1/4})$ follows from the inequality
\begin{equation*}|K|^{\frac{n-k}{n}}\ls (cL_K)^k\max_{H\in {\rm Gr}_{n-k}}|K\cap H|,\end{equation*}
where $L_K$ is the isotropic constant of $K$ (see e.g. \cite[Proposition~5.1]{Chasapis-Giannopoulos-Liakopoulos-2015}),
and the estimate $L_K=O(n^{1/4})$ of Klartag \cite{Klartag-2006} who improved an earlier result of  Bourgain \cite{Bourgain-1991}.
For several special classes of bodies the isotropic constant is uniformly bounded, and hence the answer
to the slicing problem is known to be affirmative; see \cite{BGVV-book}.

\medskip

In the case where $K$ is a generalized $k$-intersection body
in $\R^n$ (we write $K\in {\cal{BP}}_k^n$; see definition in Section \ref{sections})
and $L$ is any origin-symmetric star body in $\R^n,$ inequality \eqref{main-prob1} was proved in \cite{Koldobsky-2011} for $k=1$, and in \cite{Koldobsky-Ma} for $1<k<n:$
\begin{equation} \label{initial-sect}
|K|^{\frac {n-k}n}-|L|^{\frac {n-k}n}\ls c_{n,k}^k
\max_{F\in {\rm Gr}_{n-k}} \big(|K\cap F|-|L\cap F|\big),
\end{equation}
where $c_{n,k}^k=\omega_n^{\frac {n-k}n}/\omega_{n-k},$ and $\omega_n$ is the volume of
the unit Euclidean  ball in $\R^n.$
One can check that $c_{n,k}\in (\frac 1{\sqrt{e}},1)$ for all $n,k.$

Note that in Question \ref{q-sect} we added an extra assumption
that $L\subset K,$ compared to \eqref{initial-sect}. Without extra assumptions on $K$ and $L,$
inequality \eqref{main-prob1} cannot hold with any $\gamma>0,$ as follows from counterexamples
to the Busemann-Petty
problem. The Busemann-Petty problem asks whether, for any origin-symmetric convex bodies
$K$ and $L,$ inequalities $|K\cap F|\ls |L\cap F|$ for all
$F\in {\rm Gr}_{n-k}$ necessarily imply $|K|\ls |L|.$ The answer is negative in general;
see \cite[Chapter~5]{Koldobsky-book} for details. Every counterexample provides a pair of
bodies $K$ and $L$ that contradict inequality \eqref{main-prob1}. However, if $K$ is a generalized
$k$-intersection body, the answer to the question of Busemann and Petty is affirmative,
as proved by Lutwak \cite{Lutwak-1988} for $k=1,$ and by Zhang \cite{Zhang-1996} for $k>1.$
Inequality \eqref{initial-sect} is a quantified version of this fact.

\medskip

Our first result extends \eqref{initial-sect} to arbitrary origin-symmetric star bodies.
For a star body $K$ in $\R^n$ and $1\ls k <n,$ denote by
\begin{equation}\label{ovr}d_{\rm {ovr}}(K,{\cal{BP}}_k^n) = \inf \left\{ \left( \frac {|D|}{|K|}\right)^{1/n}:
\ K\subset D,\ D\in {\cal{BP}}_k^n \right\}\end{equation}
the outer volume ratio distance from $K$ to the class of generalized $k$-intersection
bodies.

\begin{theorem} \label{vdi-volume} Let $1\ls k <n,$ and let $K$ and $L$ be origin-symmetric star bodies in $\R^n$ such that
$L\subset K.$ Then
\begin{equation} \label{vdi-volume-ineq}
|K|^{\frac {n-k}n}-|L|^{\frac {n-k}n}\ls c_{n,k}^k d_{\rm ovr}^k(K,{\cal{BP}}_k^n)
\max_{F\in {\rm Gr}_{n-k}} \big(|K\cap F|-|L\cap F|\big).
\end{equation}
\end{theorem}

By John's theorem \cite{John-1948} and the fact that ellipsoids are intersection bodies,
if $K$ is origin-symmetric and convex, then $d_{\rm ovr}(K,{\cal{BP}}_k^n)\ls \sqrt{n}.$
In fact the same is true for any convex body by K.~Ball's volume ratio estimate in \cite{Ball-1991c}.
The outer volume ratio distance was also estimated in \cite{KPZ}. If $K$ is an origin-symmetric convex
body in $\R^n,$ then
\begin{equation}\label{eq:ovr-10}d_{\rm ovr}(K,{\cal{BP}}_k^n)\ls c\sqrt{n/k}\,[\log (en/k)]^{\frac{3}{2}},\end{equation}
where $c>0$ is an absolute constant.  In conjunction with Theorem \ref{vdi-volume}, this estimate provides
an affirmative answer to Question \ref{q-sect} for sections of proportional dimensions.

\begin{corollary}\label{proport} Let $1\ls k <n,$ let $K$ be an origin-symmetric convex body in $\R^n,$
and let $L$ be an origin-symmetric star body in $\R^n$ such that $L\subset K.$ Then
\begin{equation} \label{vdi-volume-ineq-2}
|K|^{\frac {n-k}n}-|L|^{\frac {n-k}n}\ls C^k \left(\sqrt{n/k}\,[\log (en/k)]^{\frac{3}{2}}\right)^k
\max_{F\in {\rm Gr}_{n-k}} \big(|K\cap F|-|L\cap F|\big),
\end{equation}
where $C$ is an absolute constant.
\end{corollary}

It is also known that for several classes of origin-symmetric convex bodies the distance $d_{\rm ovr}(K,{\cal{BP}}_k^n)$
is bounded by an absolute constant. These classes include unconditional convex bodies,
duals of bodies with bounded volume ratio (see \cite{Koldobsky-2015}) and the unit balls of normed
spaces that embed in $L_p,\ -n<p<\infty$
(see \cite{Koldobsky-2016, Milman-2006, Koldobsky-Pajor}).

\smallbreak

The inequality of Theorem \ref{vdi-volume} can be extended to arbitrary measures in place of volume,
as follows. Let $f$ be a bounded non-negative measurable function on $\R^n.$ Let $\mu$ be the measure
with density $f$ so that $\mu(B)=\int_B f$  for every Borel set $B$ in $\R^n.$ Also, for every $F\in {\rm Gr}_{n-k}$ we write
$\mu(B\cap F)=\int_{B\cap F} f,$ where we integrate the restriction of $f$ to $F$ against Lebesgue measure
on $F.$

It was proved in \cite{Koldobsky-2015} that for any $1\ls k <n,$ any origin-symmmetric star body $K$
in $\R^n$ and any measure $\mu$ with even non-negative continuous density $f$ in $\R^n,$
\begin{equation} \label{slicing-measure}
\mu(K)\ls \frac n{n-k} c_{n,k}^k\ |K|^{\frac kn}\ d_{\rm ovr}^k(K,{\cal{BP}}_k^n)
\max_{F\in {\rm Gr}_{n-k}} \mu(K\cap F).
\end{equation}

Considering measures with densities supported in $K\setminus L$ in inequality \eqref{slicing-measure},
we get the following measure difference inequality.

\begin{theorem} \label{vdi-measure} Let $1\ls k <n,$ let $K$ and $L$ be origin-symmetric star bodies in $\R^n$ such that
$L\subset K,$ and let $\mu$ be a measure with even non-negative continuous density. Then
\begin{equation} \label{vdi-measure-ineq}
\mu(K)-\mu(L)\ls \frac n{n-k} c_{n,k}^k\ |K|^{\frac kn}\ d_{\rm ovr}^k(K,{\cal{BP}}_k^n)
\max_{F\in {\rm Gr}_{n-k}} \big(\mu(K\cap F)-\mu(L\cap F)\big).
\end{equation}
\end{theorem}

In Section \ref{sections} we provide an alternative proof of this result.

\smallbreak

Moreover, using an approach recently developed in \cite{Chasapis-Giannopoulos-Liakopoulos-2015}, we
prove a different version of Theorem \ref{vdi-measure}, where the symmetry and continuity assumptions
are dropped, but the body $K$ is required to be convex.

\begin{theorem} \label{vdi-measure-no2} Let $1\ls k <n,$ let $K$ be a convex body with $0\in K$ and
let $L\subseteq K$ be a Borel set in $\R^n$.
For any measure $\mu $ with a bounded measurable non-negative density, we have
\begin{equation}
\mu (K)^{n-k}-\mu (L)^{n-k}\ls \left (c_0\sqrt{n-k}\right )^{k(n-k)}|K|^{\frac{k(n-k)}{n}}\max_{F\in G_{n,n-k}}\big (\mu (K\cap F)^{n-k}-\mu (L\cap F)^{n-k}\big )
\end{equation}
where $c_0>0$ is an absolute constant.
\end{theorem}

A different kind of volume difference inequality was proved in \cite{GK}. If $K$ is any origin-symmetric star body in $\R^n$, $L$ is an
intersection body, and
$\min_{\xi\in S^{n-1}} \big(|K\cap \xi^{\perp }|-|L\cap \xi^{\perp }|\big) > 0,$
where $\xi^\bot$ is the subspace of $\R^n$ perpendicular to $\xi$, then
\begin{equation}\label{initial-sect-low}
|K|^{\frac{n-1}n} - |L|^{\frac{n-1}n} \gr c \frac{1}{\sqrt{n}M(\overline{L})}
\min_{\xi\in S^{n-1}} \big(|K\cap \xi^{\perp }|-|L\cap \xi^{\perp }|\big),
\end{equation}
where $c>0$ is an absolute constant, $\overline{L}=L/|L|^{\frac 1n},$ $M(L)=\int_{S^{n-1}}\|\theta\|_L d\sigma(\theta),$ and
$\sigma$ is the normalized Lebesgue measure on the sphere.

As shown in \cite{GM}, there exist constants $c_1,c_2>0$ such that
for any $n\in \N$ and any origin-symmetric convex body $K$ in $\R^n$ in the isotropic position,
\begin{equation}\label{eq:sep-1}\frac{1}{M(K)} \gr\ c_1 \frac{n^{1/10}L_K}{\log^{2/5}(e+n)}\gr c_2\frac{n^{1/10}}{\log^{2/5}(e+n)}.\end{equation}
Also, if $K$ is convex, has volume $1$ and is in the minimal mean width position, then we have
\begin{equation}\label{eq:sep-22}\frac{1}{M(K)} \gr\ c_3\frac{\sqrt{n}}{\log (e+n)}.\end{equation}
Inserting these estimates into \eqref{initial-sect-low} we obtain estimates independent
from the bodies.

For a star body $K$ in $\R^n$ and $1\ls k <n,$ we define
$$d_{k}(K,{\cal{BP}}_k^n) = \inf \left\{ \left(\frac {\int_{S^{n-1}}\|\theta\|_K^{-k} d\sigma (\theta )}
{\int_{S^{n-1}}\|\theta\|_D^{-k} d\sigma (\theta )}\right)^{\frac 1k}:
\ D\subset K,\ D\in {\cal{BP}}_k^n \right\}.$$
By John's theorem, if $K$ is origin-symmetric and convex, then $d_{k}(K,{\cal{BP}}_k^n) \ls \sqrt{n}.$

We prove the following generalization of \eqref{initial-sect-low}.

\begin{theorem} \label{vdi-sect-low} Let $1\ls k <n,$ and let $K$ and $L$ be origin-symmetric star bodies in $\R^n$ such that
$L\subset K.$ Then
\begin{equation}\label{sect-low}
d_{k}^k (L,{\cal{BP}}_k^n) \left(|K|^{\frac{n-k}n} - |L|^{\frac{n-k}n}\right) \gr c^k \frac{1}{(\sqrt{n}M(\overline{L}))^k}
\min_{F\in {\rm Gr}_{n-k}} \big(|K\cap F|-|L\cap F|\big),
\end{equation}
where $c>0$ is an absolute constant.
\end{theorem}

We introduce another method that gives a different generalization of \eqref{initial-sect-low}.

\begin{theorem}\label{low} Let $1\ls k <n,$ and let $K$ and $L$ be bounded Borel sets in $\R^n$ with
$L\subset K$. Then
\begin{equation}\big (|K|-|L|\big )^{\frac{n-k}{n}}\gr c_{n,k}^k\min_{F\in {\rm Gr}_{n-k}}\,\big(|K\cap F|-|L\cap F|\big),\end{equation}
where $c_{n,k}^k=\omega_n^{\frac {n-k}n}/\omega_{n-k}$.
\end{theorem}

Note that Theorem \ref{low} holds true for an arbitrary pair of bounded Borel sets $L\subseteq K$ and
it no longer involves the distance $d_k$ and $M(\overline{L})$. Actually, the constant $c_{n,k}$ is sharp as one can check from the example
of the ball $K=B_2^n$ and $L=\beta B_2^n$ where $\beta\to 0$. Nevertheless, it is formally not stronger than Theorem \ref{vdi-sect-low} because
$|K|^{\frac{n-k}n} - |L|^{\frac{n-k}n}$ is smaller than $\left (|K|-|L|\right )^{\frac{n-k}{n}}.$

We deduce Theorem \ref{low} from a more general statement for arbitrary measures.

\begin{theorem}\label{low-measure} Let $1\ls k <n,$ and let $K$ and $L$ be two bounded Borel sets in $\R^n$ such that
$L\subset K$. Let $\mu $ a measure in ${\mathbb R}^n$ with bounded density $g$. Then,
\begin{equation}\big (\mu (K)-\mu (L)\big )^{\frac{n-k}{n}} \gr c_{n,k}^k\,\frac{1}{\|g\|_{\infty }^{\frac{k}{n}}}
\left (\int_{{\rm Gr}_{n-k}}\big (\mu (K\cap F)-\mu (L\cap F)\big )^{\frac{n}{n-k}}\,d\nu_{n,n-k}(F)\right )^{\frac{n-k}{n}},\end{equation}
where $\nu_{n,n-k}$ is the Haar probability measure on ${\rm Gr}_{n-k}$. In particular,
\begin{equation}\big (\mu (K)-\mu (L)\big )^{\frac{n-k}{n}}
\gr c_{n,k}^k\frac{1}{\|g\|_{\infty }^{\frac{k}{n}}}\min_{F\in {\rm Gr}_{n-k}}\big (\mu (K\cap F)-\mu (L\cap F)\big ).\end{equation}
\end{theorem}

\medskip

An inequality going in the direction opposite to \eqref{sect-low} was proved in \cite{Koldobsky-2015}.
Suppose that $K$ is an infinitely smooth origin-symmetric convex body in $\R^n$, with strictly positive curvature,
that is not an intersection body. Then there exists an origin-symmetric convex body $L$ in $\R^n$ such that $L\subset K$
and
\begin{equation}\label{initial-sect-counter-1}
|K|^{\frac {n-1}n}- |L|^{\frac {n-1}n} < c_{n,1}\min_{\xi\in S^{n-1}}\, \big(|K\cap \xi^\bot| - |L\cap \xi^\bot|\big ).
\end{equation}
Here we prove a similar inequality going in the direction opposite to \eqref{vdi-volume-ineq}.

\begin{theorem} \label{counter-sect-max} Suppose that $L$ is an infinitely smooth origin-symmetric
convex body in $\R^n$, with strictly positive curvature, that is not an intersection body. Then
there exists an origin-symmetric convex body $K$ in $\R^n$ such that $L\subset K$ and
\begin{equation}\label{initial-sect-counter-2}
 |K|^{\frac {n-1}n}- |L|^{\frac {n-1}n} > c\frac 1{\sqrt{n}M(\overline{L})} \max_{\xi\in S^{n-1}}
\,\big(|K\cap \xi^\bot| - |L\cap \xi^\bot|\big),
\end{equation}
where $c>0$ is an absolute constant.
\end{theorem}

\medskip

Let us pass to projections. For $\xi\in S^{n-1}$ and a convex body $L,$
we denote by $L\vert\xi^\bot$ the orthogonal projection of
$L$ to $\xi^\bot.$ Let $\beta_n$ be the smallest constant $\beta>0$ satisfying
\begin{equation}\label{question-proj}
\beta \big (|L|^{\frac {n-1}n}-|K|^{\frac {n-1}n}\big ) \gr
\min_{\xi\in S^{n-1}} \,\big(|L\vert \xi^\bot|-|K\vert \xi^\bot|\big)
\end{equation}
for all origin-symmetric convex bodies $K,L$ in $\R^n$ whose curvature
functions $f_K$ and $f_L$ exist and satisfy
$f_K(\xi)\ls f_L(\xi)$ for all $\xi\in S^{n-1}.$ We prove
\begin{theorem}\label{proj-sqrtn}
$\beta_n\simeq \sqrt{n},$
i.e. there exist absolute constants $a,b>0$ such that for all $n\in \N$
$$a\sqrt{n}\ls \beta_n \ls b\sqrt{n}.$$
\end{theorem}

It was proved in \cite{Koldobsky-2011, Koldobsky-2013} that if $L$ is a projection body (see definition
in Section \ref{projections})
and $K$ is an origin-symmetric convex body, then
\begin{equation} \label{initial-proj}
|L|^{\frac {n-1}n}-|K|^{\frac {n-1}n} \gr c_{n,1}
\min_{\xi\in S^{n-1}} \,\big(|L\vert \xi^\bot|-|K\vert \xi^\bot|\big).
\end{equation}
Note that we formulate \eqref{question-proj} with the condition $f_K\ls f_L,$ which is not needed
for \eqref{initial-proj}. The reason is that without an extra condition inequality \eqref{question-proj} simply
cannot hold in general with any $\beta>0$. This follows from counterexamples to the Shephard problem
asking whether, for any origin-symmetric convex bodies
$K$ and $L,$ inequalities $|K\vert\xi^\bot|\ls |L\vert\xi^\bot|$ for all
$\xi\in S^{n-1}$ necessarily imply $|K|\ls |L|.$ The answer is negative in general;
see \cite{Petty-1965, Schneider-1967} or \cite[Chapter~8]{Koldobsky-book} for details. However, if $L$ is a projection body,
the answer to the question of Shephard is affirmative, as proved by Petty \cite{Petty-1965} and Schneider \cite{Schneider-1967}.
Inequality \eqref{initial-proj} is a quantified version of this fact.

\medskip

For a convex body $L$ in $\R^n$ denote by
$$d_{\rm {vr}}(L,\Pi) = \inf \left\{ \left( \frac {|L|}{|D|}\right)^{1/n}:
\ D\subset L,\ D\in \Pi \right\}$$
the volume ratio distance from $L$ to the class of projection bodies.
We extend \eqref{initial-proj} to arbitrary origin-symmetric convex bodies, as follows.

\begin{theorem} \label{main-proj1} Suppose that $K$ and $L$ are origin-symmetric convex
bodies in $\R^n,$ and their curvature functions exist and satisfy $f_K(\xi)\ls f_L(\xi)$ for all $\xi\in S^{n-1}.$
Then
\begin{equation} \label{vdi-proj}
d_{\rm {vr}}(L,\Pi)\,\big(|L|^{\frac {n-1}n}-|K|^{\frac {n-1}n}\big ) \gr c_{n,1}
\min_{\xi\in S^{n-1}} \,\big(|L\vert \xi^\bot|-|K\vert \xi^\bot|\big).
\end{equation}
\end{theorem}

Again by K.~Ball's volume ratio estimate, for any convex body $K$ in $\R^n,$
$d_{\rm vr}(K,\Pi)\ls \sqrt{n}.$ In Section \ref{projections} we show that this distance
can be of the order $\sqrt{n},$ up to an absolute constant. The same argument is used to
deduce Theorem \ref{proj-sqrtn} from Theorem \ref{main-proj1}.

\medbreak

Denote by $h_K$ the support function, and by
$$w(K)=\int_{S^{n-1}} h_K(\xi) d\sigma(\xi)$$
the mean width of the body $K.$ Denote by
$$d_w(K,\Pi) = \inf\left\{\frac{w(D)}{w(K)}: K\subset D,\ D\in \Pi \right\}$$
the mean width distance from $K$ to the class of projection bodies.

\begin{theorem} \label{main-proj2} Suppose that $K$ and $L$ are origin-symmetric convex
bodies in $\R^n,$ and their curvature functions exist and satisfy $f_K(\xi)\ls f_L(\xi)$ for all $\xi\in S^{n-1}.$
Then
\begin{equation} \label{vdi-proj2}
|L|^{\frac {n-1}n}-|K|^{\frac {n-1}n} \ls c\ d_{\rm w}(K,\Pi) \frac{w(\overline{K})}{\sqrt{n}}
\max_{\xi\in S^{n-1}} \,\big(|L\vert \xi^\bot|-|K\vert \xi^\bot|\big),
\end{equation}
where $c$ is an absolute constant.
\end{theorem}
In Section \ref{projections} we show that the distance $d_w$ can be of the order $\sqrt{n},$
up to a logarithmic term. Note that if $K$ is a symmetric convex body of volume 1 in $\R^n$ and is in the minimal mean width position, then $w(K)\ls c\sqrt{n}(\log n).$

\smallbreak

Theorems \ref{main-proj1} and \ref{main-proj2} are complemented by the following results, going in the opposite directions, that were proved in \cite{Koldobsky-2016+}. The constant in Theorem \ref{counter-min} is written
in a more general form than in \cite{Koldobsky-2016+}.

\begin{theorem}\label{counter-max}
Suppose that $L$ is an origin-symmetric convex body in $\R^n$, with strictly positive curvature, 
that is not a projection body. Then there exists an origin-symmetric convex body $K$
in $\R^n$ so that $f_L(\xi)\gr f_K(\xi)$ for all $\xi\in S^{n-1}$ and
$$\max_{\xi\in S^{n-1}} \,\big(|L\vert \xi^\bot|-|K\vert \xi^\bot|\big )\ls \frac 1{c_{n,1}} \big(|L|^{\frac {n-1}n}-|K|^{\frac {n-1}n}\big).$$
\end{theorem}

\begin{theorem} \label{counter-min} Suppose that $K$ is an origin-symmetric convex body in $\R^n$
that is not a projection body.
Then there exists an origin-symmetric convex body $L$ in $\R^n$ so that $f_L(\xi)\gr f_K(\xi)$ for
all $\xi\in S^{n-1}$ and
$$\min_{\xi\in S^{n-1}}\,\big (|L\vert \xi^\bot|-|K\vert \xi^\bot|\big ) \gr \frac {c\sqrt{n}}{w(\overline{K})}\,\big (|L|^{\frac {n-1}n}-|K|^{\frac {n-1}n}\big ),$$
where $c$ is an absolute constant.
\end{theorem}

In Section \ref{sections} we provide the proofs of the volume difference inequalities for sections, and in Section \ref{projections}
we give the proofs of the volume difference inequalities for projections. As we proceed, we introduce notation and
the necessary background information. We refer to the books \cite{Gardner-book} and \cite{Schneider-book}
for basic facts from the Brunn-Minkowski theory and to the book \cite{AGA-book} for basic facts from asymptotic convex geometry.

%%%%%%%%%%%%%%%%%%%%%%%%%%%%%%%%%%%%%%%%%%%%%%%%%%%%%%%%%%%%%%%%%%%%%%%%%%%%%%%%%%%%%%%%%%%%%%%%%%%%%%%%%%%%%%%%%%%%%%%%%%%%%%%%%%%%%%%%%%
\section{Volume difference inequalities for sections}\label{sections}
%%%%%%%%%%%%%%%%%%%%%%%%%%%%%%%%%%%%%%%%%%%%%%%%%%%%%%%%%%%%%%%%%%%%%%%%%%%%%%%%%%%%%%%%%%%%%%%%%%%%%%%%%%%%%%%%%%%%%%%%%%%%%%%%%%%%%%%%%%

We need several definitions from convex geometry. A closed bounded set $K$ in $\R^n$ is called a star body if
every straight line passing through the origin crosses the boundary of $K$ at exactly two points different from the origin,
the origin is an interior point of $K$, and the Minkowski functional of $K$ defined by
\begin{equation}\label{eq:sep-2}\|x\|_K = \min\{a\gr 0:\ x\in aK\}\end{equation}
is a continuous function on $\R^n$.

The radial function of a star body $K$ is defined by
\begin{equation}\label{eq:sep-3}\rho_K(x) = \|x\|_K^{-1}, \qquad x\in \R^n,\ x\neq 0.\end{equation}
If $x\in S^{n-1}$ then $\rho_K(x)$ is the radius of $K$ in the direction of $x$.

We use the polar formula for the volume of a star body:
\begin{equation}\label{polar-volume}|K|=\frac{1}{n}\int_{S^{n-1}} \|\theta\|_K^{-n} d\theta,
\end{equation}
where $d\theta$ stands for the uniform measure on the sphere with density 1.

The class ${\cal{BP}}_k^n$ of {\it generalized $k$-intersection bodies} was introduced by Lutwak \cite{Lutwak-1988} for $k=1,$ and by Zhang \cite{Zhang-1996} for $k>1.$
For $1\ls k \ls n-1,$  the {\it $(n-k)$-dimensional spherical Radon transform}
$R_{n-k}:C(S^{n-1})\to C({\rm Gr}_{n-k})$
is a linear operator defined by
\begin{equation}\label{radon}R_{n-k}g (E)=\int_{S^{n-1}\cap E} g(\theta)\ d\theta,\qquad E\in {\rm Gr}_{n-k}\end{equation}
for every function $g\in C(S^{n-1}).$ We say that
an origin-symmetric star body $D$ in $\R^n$ is a {\it generalized $k$-intersection body},
and write $D\in {\cal{BP}}_k^n$, if there exists a finite non-negative Borel measure $\mu_D$
on ${\rm Gr}_{n-k}$ so that for every $g\in C(S^{n-1})$
\begin{equation}\label{def-bpk}
\int_{S^{n-1}} \rho_D^{k}(\theta) g(\theta)\ d\theta=\int_{{\rm Gr}_{n-k}} R_{n-k}g(H)\ d\mu_D(H).
\end{equation}
The class ${\cal{BP}}_1^n$ is the original class of intersection bodies introduced by Lutwak.

\bigskip

\noindent{\bf Proof of Theorem \ref{vdi-volume}.} For every $H\in {\rm Gr}_{n-k}$ we have
$$|K\cap H|-|L\cap H|\ls \max_{F\in {\rm Gr}_{n-k}} \left(|K\cap F|-|L\cap F|\right).$$
Writing volume in terms of the Radon transform, we get
$$\frac 1{n-k} \left(R_{n-k}(\|\cdot\|_K^{-n+k})(H) - R_{n-k}(\|\cdot\|_L^{-n+k})(H)\right)\ls
\max_{F\in {\rm Gr}_{n-k}} \left(|K\cap F|-|L\cap F|\right).$$
Let $D\in {\cal{BP}}_k^n,\ K\subset D.$ Integrating both sides by $H\in {\rm Gr}_{n-k}$ with the measure $\mu_D$
corresponding to $D$ by \eqref{def-bpk}, we get
\begin{equation} \label{intrD}
\frac 1{n-k} \int_{S^{n-1}} \|\theta\|_D^{-k}\left(\|\theta\|_K^{-n+k}-\|\theta\|_L^{-n+k}\right) d\theta\ls
\max_{F\in {\rm Gr}_{n-k}} \left(|K\cap F|-|L\cap F|\right) \mu_D({\rm Gr}_{n-k}).
\end{equation}
We have $\|\theta\|_D^{-1}\gr \|\theta\|_K^{-1} \gr \|\theta\|_L^{-1},$ because $L\subset K \subset D.$
Using this, H\"older's inequality and the polar formula for volume, we estimate the left-hand side of \eqref{intrD}
by
$$\frac 1{n-k} \int_{S^{n-1}} \|\theta\|_K^{-k}\left(\|\theta\|_K^{-n+k}-\|\theta\|_L^{-n+k}\right) d\theta
\gr \frac n{n-k}\left(|K|-|K|^{\frac kn}|L|^{\frac {n-k}n}\right).$$

To estimate $\mu_D({\rm Gr}_{n-k})$ from above, we combine the fact that $1= R_{n-k}{\bf 1}(E)/|S^{n-k-1}|$ for every $E\in {\rm Gr}_{n-k}$ with
\eqref{def-bpk} and H\"older's inequality to write
\begin{align}\label{eq:ovr-5} \mu_D({\rm Gr}_{n-k}) &= \frac 1{\left|S^{n-k-1}\right|} \int_{{\rm Gr}_{n-k}} R_{n-k}{\bf 1}(E) d\mu_D(E)\\
\nonumber &=\frac 1{\left| S^{n-k-1} \right| } \int_{S^{n-1}} \|\theta\|_D^{-k}\ d\theta \\
\nonumber &\ls  \frac 1{\left|S^{n-k-1}\right|} \left|S^{n-1}\right|^{\frac{n-k}n} \left(\int_{S^{n-1}} \|\theta\|_D^{-n}\ d\theta\right)^{\frac kn}\\
\nonumber &=  \frac{1}{\left|S^{n-k-1}\right|} \left|S^{n-1}\right|^{\frac{n-k}n} n^{\frac{k}{n}}|D|^{\frac{k}{n}}.
\end{align}
These estimates show that
\begin{align}\label{eq:ovr-6}\frac n{n-k}\left(|K|-|K|^{\frac kn}|L|^{\frac {n-k}n}\right) &\ls \frac{1}{\left|S^{n-k-1}\right|} \left|S^{n-1}\right|^{\frac{n-k}n} n^{\frac{k}{n}}|D|^{\frac{k}{n}}
\max_{F\in {\rm Gr}_{n-k}} \left(|K\cap F|-|L\cap F|\right) \\
\nonumber & =\frac{n}{n-k} \,c_{n,k}^k|D|^{\frac{k}{n}}
\max_{F\in {\rm Gr}_{n-k}} \left(|K\cap F|-|L\cap F|\right).
\end{align}
Finally, we choose $D$ so that $|D|^{1/n}\ls (1+\delta)d_{\rm ovr}(K,{\cal{BP}}_k^n)|K|^{1/n},$
and then send $\delta$ to zero. \qed

\bigskip

Next, we extend Theorem \ref{vdi-volume} to arbitrary measures in place of volume. Let $f$ be a bounded non-negative measurable function on $\R^n$
and let $\mu$ be the measure with density $f$. Writing integrals in polar coordinates, we get
\begin{equation}\label{measure-polar}
\mu(K)= \int_K f(x) dx = \int_{S^{n-1}} \left(\int_0^{\rho_K(\theta)} r^{n-1}f(r\theta) dr \right) d\theta,
\end{equation}
and for $H\in {\rm Gr}_{n-k}$
\begin{align}\label{measure-polar-lowdim}
\mu(K\cap H)&= \int_{K\cap H} f(x) dx = \int_{S^{n-1}\cap H}
\left(\int_0^{\rho_K(\theta)} r^{n-k-1}f(r\theta) dr \right) d\theta\\
\nonumber &=R_{n-k}\left(\int_0^{\rho_K(\cdot)} r^{n-k-1}f(r\cdot) dr\right)(H).
\end{align}

\bigskip

\noindent{\bf Proof of Theorem \ref{vdi-measure}.}  Let $f$ be the density of the measure $\mu.$
For every $H\in {\rm Gr}_{n-k}$ we have
$$\mu(K\cap H)-\mu(L\cap H)\ls \max_{F\in {\rm Gr}_{n-k}} \left(\mu(K\cap F)-\mu(L\cap F)\right).$$
Using \eqref{measure-polar-lowdim}, we get
$$R_{n-k}\left(\int_{\rho_L(\cdot)}^{\rho_K(\cdot)} r^{n-k-1}f(r\cdot) dr\right)(H)
\ls \max_{F\in {\rm Gr}_{n-k}} \left(\mu(K\cap F)-\mu(L\cap F)\right).$$
Let $D\in {\cal{BP}}_k^n,\ K\subset D.$ Integrating both sides by $H\in {\rm Gr}_{n-k}$ with the measure $\mu_D$
corresponding to $D$ by \eqref{def-bpk}, we get
\begin{equation} \label{intrM}
\int_{S^{n-1}} \rho_D^k(\theta)\left(\int_{\rho_L(\theta)}^{\rho_K(\theta)} r^{n-k-1}f(r\theta) dr\right)d\theta
\ls \max_{F\in {\rm Gr}_{n-k}} \left(\mu(K\cap F)-\mu(L\cap F)\right) \mu_D({\rm Gr}_{n-k}).
\end{equation}
We have $\rho_D\gr \rho_K \gr \rho_L,$ because $L\subset K \subset D.$
Using this and \eqref{measure-polar}, we estimate the left-hand side of \eqref{intrM}
from below
\begin{align*}
\int_{S^{n-1}} \rho_D^k(\theta)\left(\int_{\rho_L(\theta)}^{\rho_K(\theta)} r^{n-k-1}f(r\theta) dr\right)d\theta
&\gr \int_{S^{n-1}} \rho_K^k(\theta)\left(\int_{\rho_L(\theta)}^{\rho_K(\theta)} r^{n-k-1}f(r\theta) dr\right)d\theta\\
&\gr \int_{S^{n-1}} \left(\int_{\rho_L(\theta)}^{\rho_K(\theta)} r^{n-1}f(r\theta) dr\right)d\theta = \mu(K)-\mu(L).
\end{align*}
Now estimate $\mu_D(G_{n-k})$ and then choose $D$ in the same way as in the proof
of Theorem \ref{vdi-volume}.\qed

\begin{remark}\rm Note that in the case of volume ($f\equiv 1$), Theorem \ref{vdi-measure} implies that
if $K$ is an origin-symmetric convex body in $\R^n,$
and $L$ is an origin-symmetric star body in $\R^n$ such that $L\subset K$ then
\begin{equation*}
|K|^{\frac{n-k}n}-|L|^{\frac{n-k}n}\ls\frac{|K|-|L|}{|K|^{\frac kn}}\ls \frac n{n-k} c_{n,k}^k\,d_{\rm ovr}^k(K,{\cal{BP}}_k^n)
\max_{F\in {\rm Gr}_{n-k}} \big(|K\cap F|-|L\cap F|\big).
\end{equation*}
This estimate differs from the one of Theorem \ref{vdi-volume} by a factor $\frac{n}{n-k}$; however, note that also
$(|K|-|L|)/|K|^{\frac kn}$ is greater than $|K|^{\frac{n-k}n}-|L|^{\frac{n-k}n}.$
\end{remark}

To prove Theorem \ref{vdi-measure-no2} we use a technique that was introduced in \cite{Chasapis-Giannopoulos-Liakopoulos-2015}.
It is based on the following generalized Blaschke-Petkantschin formula (see \cite{Gardner-2007}).

\begin{lemma}\label{lem:gardner-07}Let $1\ls q\ls s\ls n$. There exists a constant $p(n,s,q)>0$ such that, for every non-negative
bounded Borel measurable function $f:({\mathbb R}^n)^q\to {\mathbb R}$,
\begin{align}\label{eq:tools-1}&\int_{{\mathbb R}^n}\cdots \int_{{\mathbb R}^n}f(x_1,\ldots ,x_q)dx_1\cdots dx_q\\
\nonumber &\hspace*{1cm} =p(n,s,q)\int_{G_{n,s}}\int_F\cdots \int_Ff(x_1,\ldots ,x_q)\,|{\rm conv}(0,x_1,\ldots ,x_q)|^{n-s}dx_1\ldots dx_q\,d\nu_{n,s}(F),
\end{align}
where $\nu_{n,s}$ is the Haar probability measure on ${\rm Gr}_s$. The exact value of the constant $p(n,s,q)$ is
\begin{equation}\label{eq:tools-2}p(n,s,q)=(q!)^{n-s}\frac{(n\omega_n)\cdots ((n-q+1)\omega_{n-q+1})}{(s\omega_s)\cdots ((s-q+1)\omega_{s-q+1})}.\end{equation}
\end{lemma}

We will also use Grinberg's inequality:  If $D$ is a bounded Borel set of positive Lebesgue measure in ${\mathbb R}^n$ then, for any $1\ls k\ls n-1$,
\begin{equation}\label{eq:main-3}\tilde{R}_k(D):=\frac{1}{|D|^{n-k}}\int_{G_{n,n-k}}|D\cap F|^n\,d\nu_{n,n-k}(F)\ls \frac{1}{|B_2^n|^{n-k}}\int_{G_{n,n-k}}|B_2^n\cap F|^n\,d\nu_{n,n-k}(F).\end{equation}
This fact was proved by Grinberg in \cite{Grinberg-1990}. It is stated for convex bodies $D$ but the proof applies to bounded Borel sets (see also \cite{Gardner-2007}).
For the Euclidean ball we have
\begin{equation}\label{eq:main-4}\tilde{R}_k(B_2^n):=\frac{1}{|B_2^n|^{n-k}}\int_{G_{n,n-k}}|B_2^n\cap F|^n\,d\nu_{n,n-k}(F)=\frac{\omega_{n-k}^n}{\omega_n^{n-k}}=c_{n,k}^{-kn},\end{equation}
where as before
\begin{equation}\label{eq:tools-14}c_{n,k}^k:=\omega_n^{\frac{n-k}{n}}/\omega_{n-k}.\end{equation}
For any $1\ls k\ls n-1$ we define
\begin{equation*}p(n,s):=p(n,s,s).\end{equation*}
It was proved in \cite{Chasapis-Giannopoulos-Liakopoulos-2015} that for every $1\ls k\ls n-1$ we have
\begin{equation}\label{eq:tools-13}[c_{n,k}^{-n}\,p(n,n-k)]^{\frac{1}{k(n-k)}}\simeq \sqrt{n-k}.\end{equation}

\smallskip

\noindent{\bf Proof of Theorem \ref{vdi-measure-no2}.} Let $g$ be the density of the measure $\mu.$
Applying Lemma \ref{lem:gardner-07} with $q=s=n-k$ for the functions $f(x_1,\ldots ,x_{n-k})=\prod_{i=1}^{n-k}g(x_i){\bf 1}_{K}(x_i)$
and $h(x_1,\ldots ,x_{n-k})=\prod_{i=1}^{n-k}g(x_i){\bf 1}_{L}(x_i)$  we get
\begin{align}\label{eq:arb-mu-1}&\mu (K)^{n-k}-\mu (L)^{n-k} = \prod_{i=1}^{n-k}\int_Kg(x_i)dx-\prod_{i=1}^{n-k}\int_Lg(x_i)dx\\
\nonumber &=p(n,n-k)\int_{G_{n,n-k}}\Big [\int_{K\cap F}\cdots\int_{K\cap F}g(x_1)\cdots g(x_{n-k})\,|{\rm conv}(0,x_1,\ldots ,x_{n-k})|^{k}dx_1\ldots dx_{n-k}\\
\nonumber &\hspace*{1.5cm}-\int_{L\cap F}\cdots\int_{L\cap F}g(x_1)\cdots g(x_{n-k})\,|{\rm conv}(0,x_1,\ldots ,x_{n-k})|^{k}dx_1\ldots dx_{n-k}\Big ]\,d\nu_{n,n-k}(F)\\
\nonumber &=p(n,n-k)\int_{G_{n,n-k}}\int_{P_{n-k}(K,L;F)}g(x_1)\cdots g(x_{n-k})\,|{\rm conv}(0,x_1,\ldots ,x_{n-k})|^{k}dx_1\ldots dx_{n-k}\,d\nu_{n,n-k}(F),
\end{align}
where
$$P_{n-k}(K,L;F)=(K\cap F)^{n-k}\setminus (L\cap F)^{n-k}.$$
Note that
$$|{\rm conv}(0,x_1,\ldots ,x_{n-k})|^{k}\ls |K\cap F|^k$$
for all $(x_1,\ldots ,x_{n-k})\in P_{n-k}(K,L;F)$ by the convexity of $K\cap F$ and the assumption that $0\in K$. Therefore,
\begin{align}
&\mu (K)^{n-k}-\mu (L)^{n-k}\\
\nonumber &\hspace*{1cm}\ls p(n,n-k)\int_{G_{n,n-k}}|K\cap F|^k\int_{P_{n-k}(K,L;F)}g(x_1)\cdots g(x_{n-k})\,dx_1\ldots dx_{n-k}\,d\nu_{n,n-k}(F)\\
\nonumber &\hspace*{1cm}= p(n,n-k)\int_{G_{n,n-k}}|K\cap F|^k[\mu (K\cap F)^{n-k}-\mu (L\cap F)^{n-k}]\,d\nu_{n,n-k}(F)\\
\nonumber &\hspace*{1cm}\ls \max_{F\in G_{n,n-k}}[\mu (K\cap F)^{n-k}-\mu (L\cap F)^{n-k}]\cdot p(n,n-k)\int_{G_{n,n-k}}|K\cap F|^k\,d\nu_{n,n-k}(F).\end{align}
From Grinberg's inequality \eqref{eq:main-3} we have
\begin{equation}\label{eq:arb-mu-2}\int_{G_{n,n-k}}|K\cap F|^k\,d\nu_{n,n-k}(F)\ls c_{n,k}^{-kn}\,|K|^{\frac{k(n-k)}{n}}.
\end{equation}
Using also \eqref{eq:tools-13} we see that
\begin{equation}\label{eq:arb-mu-6}\mu (K)^{n-k}-\mu (L)^{n-k}\ls \left (c_0\sqrt{n-k}\right )^{k(n-k)}|K|^{\frac{k(n-k)}{n}}\max_{F\in G_{n,n-k}}[\mu (K\cap F)^{n-k}-\mu (L\cap F)^{n-k}],\end{equation}
as claimed. \qed

\begin{remark}\rm Theorem \ref{vdi-measure-no2} implies \cite[Theorem~1.1]{Chasapis-Giannopoulos-Liakopoulos-2015}: 
\begin{equation}\label{cgl}
\mu (K)\ls \left (c_0\sqrt{n-k}\right )^{k}|K|^{\frac{k}{n}}\max_{F\in G_{n,n-k}}\,\mu (K\cap F)
\end{equation}
for every convex body $K$ with $0\in K$ and any measure $\mu $. Considering measures with densities 
supported in
$K\setminus L$ in (\ref{cgl}), we get the following measure difference inequality:
\begin{equation}
\mu (K)-\mu (L)\ls \left (c_0\sqrt{n-k}\right )^{k}|K|^{\frac{k}{n}}\max_{F\in G_{n,n-k}}\big (\mu (K\cap F)-\mu (L\cap F)\big )
\end{equation}
under the assumptions of Theorem \ref{vdi-measure-no2}.
\end{remark}

The next inequalities estimate the distance between volumes of two bodies in ${\mathbb R}^n$ in terms of the minimal
difference between areas of their $(n-k)$-dimensional sections.

\bigbreak

\noindent{\bf Proof of Theorem \ref{vdi-sect-low}.} For every $H\in {\rm Gr}_{n-k}$ we have
$$|K\cap H|-|L\cap H|\gr \min_{F\in {\rm Gr}_{n-k}} \left(|K\cap F|-|L\cap F|\right).$$
Writing volume in terms of the Radon transform, we get
$$\frac 1{n-k} \left(R_{n-k}(\|\cdot\|_K^{-n+k})(H) - R_{n-k}(\|\cdot\|_L^{-n+k})(H)\right)\gr
\min_{F\in {\rm Gr}_{n-k}} \left(|K\cap F|-|L\cap F|\right).$$
Let $D\in {\cal{BP}}_k^n,\ D\subset L.$ Integrating both sides by $H\in {\rm Gr}_{n-k}$ with the measure $\mu_D$
corresponding to $D$ by \eqref{def-bpk}, we get
\begin{equation} \label{intrD-2}
\frac 1{n-k} \int_{S^{n-1}} \|\theta\|_D^{-k}\left(\|\theta\|_K^{-n+k}-\|\theta\|_L^{-n+k}\right) d\theta\gr
\min_{F\in {\rm Gr}_{n-k}} \left(|K\cap F|-|L\cap F|\right) \mu_D({\rm Gr}_{n-k}).
\end{equation}
We have $\|\theta\|_D^{-1}\ls \|\theta\|_L^{-1} \ls \|\theta\|_K^{-1},$ because $D\subset L \subset K.$
Using this, H\"older's inequality and the polar formula for volume, we estimate the left-hand side of \eqref{intrD-2}
from above by
$$\frac 1{n-k} \int_{S^{n-1}} \|\theta\|_L^{-k}\left(\|\theta\|_K^{-n+k}-\|\theta\|_L^{-n+k}\right) d\theta
\ls \frac n{n-k}\left(|L|^{\frac kn}|K|^{\frac {n-k}n}-|L|\right).$$

To estimate $\mu_D({\rm Gr}_{n-k})$ from below, we combine the fact that $1= R_{n-k}{\bf 1}(E)/|S^{n-k-1}|$ 
for every $E\in {\rm Gr}_{n-k}$ with
\eqref{def-bpk} to write
\begin{equation}\label{eq:ovr-7} \mu_D({\rm Gr}_{n-k})=
\frac 1{\left|S^{n-k-1}\right|} \int_{{\rm Gr}_{n-k}} R_{n-k}{\bf 1}(E) d\mu_D(E)=\frac {|S^{n-1}|}{\left| S^{n-k-1} \right| }
 \int_{S^{n-1}} \|\theta\|_D^{-k}\ d\sigma(\theta).
\end{equation}
These estimates show that
$$ \frac n{n-k}\left(|L|^{\frac kn}|K|^{\frac {n-k}n}-|L|\right)
\gr \frac{|S^{n-1}|}{\left|S^{n-k-1}\right|} \int_{S^{n-1}} \|\theta\|_D^{-k} d\sigma(\theta)
\min_{F\in {\rm Gr}_{n-k}} \left(|K\cap F|-|L\cap F|\right).
$$
Finally, for $\delta>0,$ we choose $D$ so that
$$\int_{S^{n-1}} \|\theta\|_D^{-k} d\sigma(\theta)\gr \frac 1{(1+\delta) d_{k}^k(L,{\cal{BP}}_k^n)}
\int_{S^{n-1}} \|\theta\|_L^{-k} d\sigma(\theta),$$
and send $\delta$ to zero. Then use Jensen's inequality and homogeneity to get
\begin{equation}\left(\int_{S^{n-1}}\|\theta\|_L^{-k}\,d\sigma(\theta)\right)^{\frac 1k}\gr
\left (\int_{S^{n-1}}\|\theta\|_L d\sigma(\theta)\right )^{-1}=\frac{1}{M(\overline{L})}|L|^{\frac{1}{n}},\end{equation}
and apply standard estimates for the $\Gamma$-function. \qed

\bigbreak

Next we prove Theorem \ref{low-measure}, which directly implies Theorem \ref{low}. For the proof we will use some basic facts
about Sylvester-type functionals. Let $C$ be a bounded Borel set of positive measure
in ${\mathbb R}^m$. For every $p>0$ we consider
the normalized $p$-th moment of the expected volume of the random simplex ${\rm conv}(0,x_1,\ldots ,x_m)$, the convex hull of the origin and
$m$ points from $C$, defined by
\begin{equation}\label{eq:tools-4}S_p(C)=\left (\frac{1}{|C|^{m+p}}
\int_C\cdots\int_C|{\rm conv}(0,x_1,\ldots ,x_m)|^pdx_1\cdots
dx_m\right )^{1/p}.\end{equation}
It was proved by Pfiefer \cite{Pfiefer-1990} (see also \cite{Gardner-2007}) that
\begin{equation*}S_p(C)\gr S_p(B_2^m).\end{equation*}
More generally, for any Borel probability measure $\nu $ on ${\mathbb R}^m$, for any $1\ls q\ls m$ and every $p>0$, we define
\begin{equation}\label{eq:tools-5}S_{p,q}(\nu )=\left (\int_{{\mathbb R}^m}\cdots\int_{{\mathbb R}^m}|{\rm conv}(0,x_1,\ldots ,x_q)|^pd\nu (x_1)\cdots
d\nu (x_q)\right )^{1/p}.\end{equation}
A generalization of Pfiefer's result appears in \cite{Dann-Paouris-Pivovarov-2015}. Let $\nu $ be a measure in ${\mathbb R}^n$ with
a bounded non-negative measurable density $g$. Then
\begin{equation}\label{eq:DPP-2015}S_{p,q}^p(\nu )\gr
\frac{\|g\|_1^{q+\frac{pq}{m}}}{\omega_m^{q+\frac{pq}{m}}\|g\|_{\infty }^\frac{pq}{m}}S_{p,q}^p({\bf 1}_{B_2^m}).\end{equation}

\medskip

\noindent{\bf Proof of Theorem \ref{low-measure}.}
Let $u(x)=g(x){\bf 1}_K(x)$ and $v(x)=g(x){\bf 1}_L(x)$. Using Lemma \ref{lem:gardner-07} with $s=n-k$ and $q=1$, we start by writing
\begin{align}\label{eq:arb-mu-11}&\mu (K)-\mu (L) = \int_{{\mathbb R}^n}u(x)dx-\int_{{\mathbb R}^n}v(x)dx\\
\nonumber &=p(n,n-k,1)\int_{G_{n,n-k}}\Big [\int_{K\cap F}g(x)\,\|x\|_2^{k}dx
-\int_{L\cap F}g(x)\,\|x\|_2^{k}dx\Big ]\,d\nu_{n,n-k}(F)\\
\nonumber &= p(n,n-k,1)\int_{G_{n,n-k}}\int_{(K\cap F)\setminus (L\cap F)}g(x)\,\|x\|_2^{k}dx\,d\nu_{n,n-k}(F).
\end{align}
(Note that $|{\rm conv}(0,x)|=\|x\|_2$, the Euclidean norm of $x$).
For every $F$ set $C_F=(K\cap F)\setminus (L\cap F)$ and consider the measure $\nu_F$ with density $g$ on $C_F$. Applying \eqref{eq:DPP-2015} with $p=k$, $q=1$ and $m=n-k$
we have
\begin{align}\label{eq:arb-mu-12}&\mu (K)-\mu (L) \gr p(n,n-k,1)\int_{{\rm Gr}_{n-k}}S_{k,1}^k(\nu_F)\,d\nu_{n,n-k}(F)\\
\nonumber &\gr p(n,n-k,1)\int_{{\rm Gr}_{n-k}}\frac{\|g\,|_{C_F}\|_1^{1+\frac{k}{n-k}}}{\omega_{n-k}^{1+\frac{k}{n-k}}\|g\,|_{C_F}\|_{\infty }^{\frac{k}{n-k}}}S_k^k({\bf 1}_{B_2^{n-k}})\,d\nu_{n,n-k}(F)\\
\nonumber &= \frac{p(n,n-k,1)}{\omega_{n-k}^{\frac{n}{n-k}}}S_2^k({\bf 1}_{B_2^{n-k}})\int_{{\rm Gr}_{n-k}}\frac{\|g\,|_{C_F}\|_1^{\frac{n}{n-k}}}{\|g\,|_{C_F}\|_{\infty }^{\frac{k}{n-k}}}\,d\nu_{n,n-k}(F).
\end{align}
Note that
\begin{equation*}p(n,n-k,1)=\frac{n\omega_n}{(n-k)\omega_{n-k}}\end{equation*}
and
\begin{equation*}S_{k,1}^k({\bf 1}_{B_2^{n-k}})=\int_{B_2^{n-k}}\|x\|_2^kdx=\frac{n-k}{n}\omega_{n-k}.\end{equation*}
Therefore,
\begin{equation*} \frac{p(n,n-k,1)}{\omega_{n-k}^{\frac{n}{n-k}}}S_2^k({\bf 1}_{B_2^{n-k}})=\frac{\omega_n}{\omega_{n-k}^{\frac{n}{n-k}}}=c_{n,k}^{\frac{kn}{n-k}}.
\end{equation*}
On the other hand, for any $F\in {\rm Gr}_{n-k}$ we have
\begin{equation*}\|g\,|_{C_F}\|_1=\mu (K\cap F)-\mu (L\cap F)\end{equation*}
and
\begin{equation*}\|g\,|_{C_F}\|_{\infty }\ls \|g\|_{\infty }.\end{equation*}
Combining the above we get
\begin{equation*}\mu (K)-\mu (L) \gr c_{n,k}^{\frac{kn}{n-k}}\frac{1}{\|g\|_{\infty }^{\frac{k}{n-k}}}\int_{{\rm Gr}_{n-k}}(\mu (K\cap F)-\mu (L\cap F)^{\frac{n}{n-k}}\,d\nu_{n,n-k}(F),\end{equation*}
and the result follows. \qed

\begin{remark}\rm Theorem \ref{low} is an immediate consequence of Theorem \ref{low-measure}. It corresponds to the case $g\equiv {\bf 1}$, for
which we clearly have $\|g\|_{\infty }=1$.
\end{remark}

We pass to Theorem \ref{counter-sect-max}. We consider Schwartz distributions, i.e. continuous functionals on the space ${\cal{S}}(\R^n)$
of rapidly decreasing infinitely differentiable functions on $\R^n$.
The Fourier transform of a distribution $f$ is defined by $\langle\hat{f}, \phi\rangle= \langle f, \hat{\phi} \rangle$ for
every test function $\phi \in {\cal{S}}(\R^n).$ For any even distribution $f$, we have $(\hat{f})^\wedge
= (2\pi)^n f$.

If $K$ is an origin-symmetric convex body and $0<p<n,$
then $\|\cdot\|_K^{-p}$  is a locally integrable function on $\R^n$ and represents a distribution acting by integration.
Suppose that $K$ is infinitely smooth, i.e. $\|\cdot\|_K\in C^\infty(S^{n-1})$ is an infinitely differentiable
function on the sphere. Then by \cite[Lemma 3.16]{Koldobsky-book}, the Fourier transform of $\|\cdot\|_K^{-p}$
is an extension of some function $g\in C^\infty(S^{n-1})$ to a homogeneous function of degree
$-n+p$ on $\R^n.$ When we write $\left(\|\cdot\|_K^{-p}\right)^\wedge(\xi),$ we mean $g(\xi),\ \xi \in S^{n-1}.$

For $f\in C^\infty(S^{n-1})$ and $0<p<n$, we denote by
$$(f\cdot r^{-p})(x) = f(x/\|x\|_2) \|x\|_2^{-p}$$
the extension of $f$ to a homogeneous function of degree $-p$ on $\R^n.$
Again by  \cite[Lemma 3.16]{Koldobsky-book}, there exists $g\in C^\infty(S^{n-1})$ such that
$$(f\cdot r^{-p})^\wedge = g\cdot r^{-n+p}.$$

If $K,L$ are infinitely smooth origin-symmetric convex bodies, the following spherical version of Parseval's
formula can be found in \cite[Lemma 3.22]{Koldobsky-book}:  for any $p\in (-n,0)$
\begin{equation}\label{parseval}
\int_{S^{n-1}} \left(\|\cdot\|_K^{-p}\right)^\wedge(\xi) \left(\|\cdot\|_L^{-n+p}\right)^\wedge(\xi) =
(2\pi)^n \int_{S^{n-1}} \|x\|_K^{-p} \|x\|_L^{-n+p}\ dx.
\end{equation}

It was proved in \cite[Theorem 1]{Koldobsky-1998} that an origin-symmetric convex body $K$ in $\R^n$ is an
intersection body if and only if the function $\|\cdot\|_K^{-1}$ represents a positive definite
distribution. In the case where $K$ is infinitely smooth, this means that the function $(\|\cdot\|_K^{-1})^\wedge$
is non-negative on the sphere.

We also need a result from \cite{Koldobsky-sectformula} (see also \cite[Theorem 3.8]{Koldobsky-book})
expressing volume of central hyperplane sections
in terms of the Fourier transform. For any origin-symmetric star
body $K$ in $\R^n,$ the distribution $(\|\cdot\|_K^{-n+1})^\wedge$ is a continuous function
on the sphere extended to a homogeneous function of degree -1 on the whole of $\R^n,$
and for every $\xi\in S^{n-1},$
\begin{equation} \label{sect-Fourier}
|K\cap \xi^\bot| = \frac 1{\pi(n-1)} (\|\cdot\|_K^{-n+1})^\wedge(\xi).
\end{equation}
In particular, if $K=B_2^n$ then for every $\xi\in S^{n-1}$
\begin{equation}\label{eucl1}
(\|\cdot\|_2^{-n+1})^\wedge(\xi) = \pi(n-1)|B_2^{n-1}|.
\end{equation}

Note that every non-intersection body can be approximated in
the radial metric by infinitely smooth non-intersection bodies with strictly positive curvature; see
\cite[Lemma~4.10]{Koldobsky-book}. Different examples of convex bodies that are not intersection bodies
(in dimensions five and higher, as in dimensions up to four such examples do not exist) can
be found in \cite[Chapter~4]{Koldobsky-book}. In particular, the unit balls of the spaces $\ell_q^n,\ q>2,\ n\gr 5$
are not intersection bodies.
\medbreak

\noindent{\bf Proof of Theorem \ref{counter-sect-max}.}
Since $L$ is infinitely smooth, the Fourier transform of $\|\cdot\|_L^{-1}$
is a continuous function on the sphere $S^{n-1}.$ Also, $L$ is not an intersection body,
so $\left(\|\cdot\|_L^{-1}\right)^\wedge < 0$ on an open set $\Omega\subset S^{n-1}.$
Let $\phi\in C^\infty(S^{n-1})$ be an even non-negative, not identically zero, infinitely smooth function
on $S^{n-1}$ with support in $\Omega\cup -\Omega.$ Extend $\phi$ to an even homogeneous
 of degree -1 function $\phi\cdot r^{-1}$ on $\R^n\setminus \{0\}.$ The Fourier transform of this
function in the sense of distributions is $\psi\cdot r^{-n+1}$ where $\psi$ is an infinitely smooth
function on the sphere.

Let $\e$ be a number such that $|B_2^{n-1}| \cdot\|\theta\|_L^{-n+1}> \e>0$ for every $\theta\in S^{n-1}.$
Define a star body $K$ by
\begin{equation}\label{newbody}
\|\theta\|_K^{-n+1}= \|\theta\|_L^{-n+1} - \delta\psi(\theta) + \frac{\e}{|B_2^{n-1}|},\qquad \theta\in S^{n-1},
\end{equation}
where $\delta>0$ is small enough so that for every $\theta$
$$|\delta\psi(\theta)|<\min\left\{\|\theta\|_L^{-n+1}-\frac{\e}{|B_2^{n-1}|},\ \frac{\e}{|B_2^{n-1}|}\right\}.$$
The latter condition implies that $L\subset K.$ Since $L$ has strictly positive curvature, by an argument
from \cite[p.~96]{Koldobsky-book}, we can make $\e, \delta$ smaller (if necessary) to ensure that the body $K$ is convex.

Now we extend the functions in \eqref{newbody} from the sphere to $\R^n\setminus \{0\}$
as homogeneous functions of degree $-n+1$ and apply the Fourier transform. We get
that for every $\xi\in S^{n-1}$
\begin{equation}\label{four0}
\left(\|\cdot\|_K^{-n+1}\right)^{\wedge}(\xi) = \left(\|\cdot\|_L^{-n+1}\right)^\wedge(\xi) - (2\pi)^n \delta \phi(\xi)
+ \pi(n-1)\e.
\end{equation}
Here, we used \eqref{eucl1} to compute the last term.
By \eqref{four0}, \eqref{sect-Fourier} and the fact that the function $\phi$ is non-negative and is equal to zero at some points,
we have
\begin{equation}\label{sect1}
 \e = \max_{\xi\in S^{n-1}}(|K\cap \xi^\bot|-|L\cap \xi^\bot|).
\end{equation}
Multiplying both sides of \eqref{four0} by $\left (\|\cdot\|_L^{-1}\right)^\wedge(\xi),$
integrating over $S^{n-1}$ and using Parseval's formula on the sphere,  we get
\begin{align*}(2\pi)^n \int_{S^{n-1}} \|\theta\|_L^{-1}\ \|\theta\|_K^{-n+1} d\theta
 &=(2\pi)^n n|L| - (2\pi)^n \delta \int_{S^{n-1}} \phi(\theta) \left (\|\cdot\|_L^{-1}\right)^\wedge(\theta) d\theta\\
 &\hspace*{1cm} +\pi(n-1) \e \int_{S^{n-1}} \left (\|\cdot\|_L^{-1}\right)^\wedge(\theta) d\theta.
 \end{align*}
Since $\phi$ is a non-negative function supported in $\Omega,$ where $\left (\|\cdot\|_L^{-1}\right)^\wedge$ is negative,
the latter equality implies
\begin{align*}(2\pi)^n n|L| + \pi(n-1) \e \int_{S^{n-1}} \left (\|\cdot\|_L^{-1}\right)^\wedge(\theta) d\theta
&< (2\pi)^n \int_{S^{n-1}} \|\theta\|_L^{-1}\ \|\theta\|_K^{-n+1} d\theta\\
&\ls (2\pi)^n \left(\int_{S^{n-1}} \|\theta\|_K^{-n} d\theta \right)^{\frac{n-1}n}
\left(\int_{S^{n-1}}\|\theta\|_L^{-n}d\theta\right)^{\frac 1n}\\
&= (2\pi)^n n |L|^{\frac 1n}|K|^{\frac {n-1}{n}}.
\end{align*}
Finally, by \eqref{eucl1}, Parseval's formula and Jensen's inequality,
\begin{align*}
\pi(n-1)\int_{S^{n-1}} \left (\|\cdot\|_L^{-1}\right)^\wedge(\theta) d\theta
&=\frac 1{|B_2^{n-1}|}\int_{S^{n-1}} \left (\|\cdot\|_L^{-1}\right)^\wedge(\theta) \left(\|\cdot\|_2^{-n+1}\right)^\wedge(\theta) d\theta\\
&=\frac {(2\pi)^n |S^{n-1}|}{|B_2^{n-1}|}\int_{S^{n-1}} \|\theta\|_L^{-1}d\sigma(\theta)\\
&\gr \frac {(2\pi)^n |S^{n-1}|}{|B_2^{n-1}|} \frac{1}{M(\overline{L})}|L|^{\frac{1}{n}}\\
&\gr c\frac {(2\pi)^n\sqrt{n} |L|^{\frac{1}{n}}}{M(\overline{L})}.
\end{align*}
Combining these estimates we get
$$(2\pi)^n n|L| + c\e\frac {(2\pi)^n\sqrt{n} |L|^{\frac{1}{n}}}{M(\overline{L})}\ls (2\pi)^n n |L|^{\frac 1n}|K|^{\frac {n-1}{n}}.$$
The result follows after we recall \eqref{sect1}. \qed

%%%%%%%%%%%%%%%%%%%%%%%%%%%%%%%%%%%%%%%%%%%%%%%%%%%%%%%%%%%%%%%%%%%%%%%%%%%%%%%%%%%%%%%%%%%%%%%%%%%%%%%%%%%%%%%%%%%%%%%%%%%%%
\section{Volume difference inequalities for projections}\label{projections}
%%%%%%%%%%%%%%%%%%%%%%%%%%%%%%%%%%%%%%%%%%%%%%%%%%%%%%%%%%%%%%%%%%%%%%%%%%%%%%%%%%%%%%%%%%%%%%%%%%%%%%%%%%%%%%%%%%%%%%%%%%%%%

The {\it support function} of a convex body $K$ in $\R^n$ is defined by
$$h_K(x) = \max\{ \langle x,y\rangle :y\in K\},\qquad x\in \R^n.$$
If $K$ is origin-symmetric, then $h_K$ is a norm on $\R^n.$

The {\it surface area measure} $S(K, \cdot)$ of a convex body $K$ in
$\R^n$ is defined as follows. For every Borel set $E \subset S^{n-1},$
$S(K,E)$ is equal to Lebesgue measure of the part of the boundary of $K$
where normal vectors belong to $E.$
We usually consider bodies with absolutely continuous surface area measures.
A convex body $K$ is said to have the {\it curvature function}
$$ f_K: S^{n-1} \to \R,$$
if its surface area measure $S(K, \cdot)$ is absolutely
continuous with respect to Lebesgue measure $\sigma_{n-1}$ on
$S^{n-1}$, and
$$
\frac{d S(K, \cdot)}{d \sigma_{n-1}}=f_K \in L_1(S^{n-1}),
$$
so $f_K$ is the density of $S(K,\cdot).$

By the approximation argument of \cite[Theorem~3.3.1]{Schneider-book},
we may assume in the formulation of Shephard's problem that the bodies
$K$ and $L$ are such that  their support functions $h_K,\,h_L$ are
infinitely smooth functions on $\R^n\setminus \{0\}$.
Using \cite[Lemma~3.16]{Koldobsky-book}
we get in this case that
the Fourier transforms $\widehat{h_K},\ \widehat{h_L}$ are the
extensions of infinitely differentiable functions on the sphere
to homogeneous distributions on $\R^n$ of degree $-n-1.$
Moreover, by a similar approximation argument (see e.g. \cite[Section 5]{GZ}),
we may assume that  our bodies have absolutely continuous surface area
measures. Therefore, in the rest of this section, $K$ and $L$ are
convex symmetric bodies with infinitely smooth support functions and absolutely
continuous surface area measures.

The following version of Parseval's formula was proved in \cite{KRZ} (see also \cite[Lemma 8.8]{Koldobsky-book}):
\begin{equation} \label{pars-proj}
\int_{S^{n-1}} \widehat{h_K} (\xi) \widehat{f_L}(\xi)\ d\xi =
(2\pi)^n \int_{S^{n-1}} h_K(x) f_L(x)\ dx.
\end{equation}

The volume of a body can be expressed in terms of its support function and
curvature function:
\begin{equation}\label{vol-proj}
|K| = \frac 1n \int_{S^{n-1}}h_K(x) f_K(x)\ dx.
\end{equation}

If $K$ and $L$ are two convex bodies in $\R^n$ the {\it mixed volume} $V_1(K,L)$
is equal to
$$V_1(K,L)= \frac{1}{n} \lim_{\e\to +0}
\frac{|K+\epsilon L|- |K|}{\e}.$$
We use the following
first Minkowski inequality (see \cite{Schneider-book} or \cite[p.23]{Koldobsky-book}):
for any convex bodies $K,L$ in $\R^n,$
\begin{equation} \label{firstmink}
V_1(K,L) \gr |K|^{\frac {n-1}n} |L|^{\frac 1n}.
\end{equation}
The mixed volume $V_1(K,L)$ can also be expressed in terms of the support and
curvature functions:

\begin{equation}\label{mixvol-proj}
V_1(K,L) = \frac 1n \int_{S^{n-1}}h_L(x) f_K(x)\ dx.
\end{equation}

Let $K$ be an origin-symmetric convex body in $\R^n.$ The {\it
projection body} $\Pi K$ of $K$ is defined as an origin-symmetric convex
body in $\R^n$ whose support function in every direction is equal to
the volume of the hyperplane projection of $K$ to this direction:
for every $\xi\in S^{n-1},$
\begin{equation} \label{def:proj}
h_{\Pi K}(\xi ) = |K\vert\xi^{\perp}|.
\end{equation}
If $L$ is the projection body of some convex body, we simply say
that $L$ is a projection body.  The Minkowski (vector) sum of projection bodies
is also a projection body. Every projection body is the limit in the Hausdorff metric
of Minkowski sums of symmetric intervals. An origin-symmetric convex body
in $\R^n$ is a projection body if and only if its polar body is the unit ball of
an $n$-dimensional subspace of $L_1;$ see \cite{Schneider-book,Gardner-book,Koldobsky-book} for proofs and
more properties of projection bodies.

\bigskip

\noindent {\bf Proof of Theorem \ref{main-proj1}.} By approximation (see \cite[Theorem 3.3.1]{Schneider-book}),
we can assume that $K,L$ are infinitely smooth. We have
\begin{equation}\label{proj11}
|L\vert \xi^\bot|- |K\vert \xi^\bot|\gr \min_{\eta\in S^{n-1}} (|L\vert \eta^\bot|-|K\vert \eta^\bot|).
\end{equation}
It was proved in \cite{KRZ} that
\begin{equation} \label{f-proj}
|K\vert \xi^\bot| = -\frac 1{\pi} \widehat{f_K}(\xi),\qquad \xi\in S^{n-1},
\end{equation}
where $f_K$ is extended from the sphere to a homogeneous function of degree
$-n-1$ on the whole $\R^n.$
Therefore, \eqref{proj11} can be written as
\begin{equation} \label{fourier-proj}
-\frac 1{\pi} \widehat{f_L}(\xi) + \frac 1{\pi} \widehat{f_K}(\xi) \gr \min_{\eta\in S^{n-1}} (|L\vert \eta^\bot|-|K\vert \eta^\bot|), \qquad \xi\in S^{n-1}.
\end{equation}
Let $D$ be a projection body such that $D\subset L,$ then $h_D\ls h_L$ in every direction.
It was proved in \cite{KRZ} that an infinitely smooth origin-symmetric convex body
$D$ in $\R^n$ is a projection body if and only if
$\widehat{h_D} \ls 0$ on the sphere $S^{n-1}.$ Integrating \eqref{fourier-proj}
with respect to this negative density, we get
$$-\int_{S^{n-1}} \widehat{h_D}(\xi) \widehat{f_L}(\xi)\ d\xi + \int_{S^{n-1}} \widehat{h_D}(\xi) \widehat{f_K}(\xi)\ d\xi \ls
\pi \int_{S^{n-1}} \widehat{h_D}(\xi)d\xi\min_{\eta\in S^{n-1}} (|L\vert \eta^\bot|-|K\vert \eta^\bot|).$$
Using Parseval's formula \eqref{pars-proj}, we get
\begin{equation}\label{proj12}
(2\pi)^n \int_{S^{n-1}} h_D(\xi) (f_L(\xi)-f_K(\xi)) d\xi \gr -\pi  \int_{S^{n-1}} \widehat{h_D}(\xi)d\xi\min_{\eta\in S^{n-1}} (|L\vert \eta^\bot|-|K\vert \eta^\bot|).
\end{equation}
We estimate the left-hand side of \eqref{proj12} from above using \eqref{vol-proj} and \eqref{mixvol-proj}
(recall that $f_K\ls f_L)$:
\begin{align}\label{proj13}
(2\pi)^n \int_{S^{n-1}} h_D(\xi) (f_L(\xi)-f_K(\xi)) d\xi &\ls (2\pi)^n \int_{S^{n-1}} h_L(\xi) (f_L(\xi)-f_K(\xi)) d\xi \\
\nonumber &\ls (2\pi)^n n (|L|-|K|^{\frac{n-1}n}|L|^{\frac 1n}).
\end{align}

To estimate the right-hand side of \eqref{proj13} from below, note that, by \eqref{f-proj},
the Fourier transform of the curvature function $f_2$ of the unit Euclidean ball is equal to
$$\widehat{f_2}(\xi) = -\pi |B_2^{n-1}|,\qquad \xi\in S^{n-1}.$$
Therefore, by \eqref{pars-proj} and \eqref{mixvol-proj} (recall that $f_2\equiv 1$) ,
\begin{align*}
-\pi \int_{S^{n-1}}\widehat{h_D}(\xi)\ d\xi &=  \frac {1}{|B_2^{n-1}|}
\int_{S^{n-1}}\widehat{h_D}(\xi)\widehat{f_2}(\xi)\ d\xi = \frac {(2\pi)^n}{|B_2^{n-1}|}
\int_{S^{n-1}} h_D(x) f_2(x)\ dx\\
&= \frac {(2\pi)^n}{|B_2^{n-1}|} nV_1(B_2^n,D)\gr  \frac {(2\pi)^n n}{|B_2^{n-1}|} |D|^{\frac 1n}
 |B_2^n|^{\frac {n-1}n}\\
& = (2\pi)^n n\ c_{n,1} |D|^{\frac 1n}.
\end{align*}
Now for $\delta>0$ choose $D$ so that
$ (1+\delta)\ d_{\rm vr}(L,\Pi)\ |D|^{\frac 1n}\gr |L|^{\frac 1n}.$
Combine the resulting inequality with \eqref{proj12} and \eqref{proj13} and send $\delta$ to zero.
 \qed

\bigskip

\noindent{\bf Proof of Theorem \ref{proj-sqrtn}.} Putting $K=\delta B_2^n$ in \eqref{question-proj}
and sending $\delta$ to zero, we get
$$\beta |L|^{\frac {n-1}n} \gr
\min_{\xi\in S^{n-1}} |L\vert \xi^\bot|.$$
By a result of K.~Ball \cite{Ball-1991},  there exists an absolute constant $c_1$ so that for each $n\in \N$
there is an origin-symmetric convex body $L_n$ in $\R^n$ satisfying
$$\min_{\xi\in S^{n-1}} |L_n\vert\xi^\bot| \gr c_1\sqrt{n} |L_n|^{\frac{n-1}n}.$$
This shows that $\beta_n\gr c_1\sqrt{n}.$
On the other hand, since ellipsoids are projection bodies, we have 
$d_{\rm vr}(L,\Pi)\ls \sqrt{n}$ for every origin-symmetric convex body $L$ in $\R^n.$
By approximation (see \cite{GZ}), one can assume that each of the bodies $L_n$
has a curvature function, so we can apply Theorem \ref{main-proj1} to the
bodies $L_n$ and $K=\delta B_2^n,\ \delta \to 0,$ to see that $\beta_n\ls (1/c_{n,1})\sqrt{n}< \sqrt{en}.$ \qed

\begin{remark}\rm From Theorem \ref{main-proj1} we see that the bodies $L_n$ defined in the proof
of Theorem \ref{proj-sqrtn} satisfy
\begin{equation*}
d_{\rm {vr}}(L_n,\Pi)|L_n|^{\frac {n-1}n}\gr c_{n,1}
\min_{\xi\in S^{n-1}}|L\vert \xi^\bot|\gr c_{n,1}\,c_1\sqrt{n} |L_n|^{\frac{n-1}n}.
\end{equation*}
This shows that $d_{\rm vr}(L_n,\Pi)\gr c_1\sqrt{n/e}$, and hence
\begin{equation*}\sup_L\,d_{{\rm vr}}(L,\Pi_n)\simeq \sqrt{n},\end{equation*}
where the supremum is over all origin-symmetric convex bodies $L$ in ${\mathbb R}^n$.
\end{remark}

\noindent {\bf Proof of Theorem \ref{main-proj2}.} Again, by approximation,
we can assume that $K,L$ are infinitely smooth. Let $D$ be a projection body such that $K\subset D,$
then $h_K\ls h_D$ in every direction. Similarly to the proof of Theorem \ref{main-proj1},
\begin{equation}\label{proj21}
(2\pi)^n \int_{S^{n-1}} h_D(\xi) (f_L(\xi)-f_K(\xi)) d\xi \ls -\pi  \int_{S^{n-1}} \widehat{h_D}(\xi)d\xi\max_{\eta\in S^{n-1}} (|L\vert \eta^\bot|-|K\vert \eta^\bot|).
\end{equation}
We estimate the left-hand side of \eqref{proj21} from below using \eqref{vol-proj} and \eqref{mixvol-proj}
(recall that $f_K\ls f_L$ and $h_K\ls h_D)$:
\begin{align}\label{proj22}
(2\pi)^n \int_{S^{n-1}} h_D(\xi) (f_L(\xi)-f_K(\xi)) d\xi &\gr (2\pi)^n \int_{S^{n-1}} h_K(\xi) (f_L(\xi)-f_K(\xi)) d\xi \\
\nonumber &\gr (2\pi)^n n (|L|^{\frac{n-1}n}|K|^{\frac 1n}-|K|).
\end{align}
Now for $\delta>0$ choose $D$ so that $$w(D)\ls (1+\delta)d_{\rm w}(K,\Pi) w(\overline{K})|K|^{\frac 1n}.$$
As in the proof of Theorem \ref{main-proj1},
\begin{align*}-\pi \int_{S^{n-1}}\widehat{h_D}(\xi)\ d\xi &= \frac {(2\pi)^n}{|B_2^{n-1}|}
\int_{S^{n-1}} h_D(x)\ dx = \frac {(2\pi)^n|S^{n-1}|}{|B_2^{n-1}|}w(D)\\
&\ls  (1+\delta) (2\pi)^nc\ d_{\rm w}(K,\Pi)\sqrt{n}
\ w(\overline{K})|K|^{\frac 1n}.
\end{align*}
We get the result combining the latter with \eqref{proj21} and \eqref{proj22} and sending $\delta$ to zero. \qed

\medskip

Finally, we show that the distance $d_w$ can be of the order $\sqrt{n},$ up to a logarithmic term.
We will use the fact that projection bodies have positions with
``small diameter". More precisely, we have the following statement: For every $D\in \Pi $ there exists $T\in GL(n)$ such that
\begin{equation}\label{eq:radius-projection-body}R(T(D))\ls\frac{\sqrt{n}}{2}|T(D)|^{1/n}.\end{equation}
In particular, this holds true if $T$ is chosen so that $T(D)$ in Lewis or L\"{o}wner or minimal mean width
position (see e.g. \cite[Chapter~4]{BGVV-book}). Let $K=B_1^n$ be the cross-polytope, and consider a projection
body $D$ such that $B_1^n\subseteq D$. We may find $T$ so that \eqref{eq:radius-projection-body} is satisfied.
We will use the next well-known result of B\'{a}r\'{a}ny and F\"{u}redi from \cite{Barany-Furedi-1988}: if $x_1,\ldots ,x_N\in RB_2^n$
then
\begin{equation*}|{\rm conv}\{x_1,\ldots ,x_N\}|^{1/n}\ls \frac{c_3R\sqrt{\log (1+N/n)}}{n}.\end{equation*}
Since
\begin{equation*}T(B_1^n)={\rm conv}\{\pm Te_1,\ldots ,\pm Te_n\}\subseteq R(T(D))B_2^n,\end{equation*}
we get
\begin{equation*}|T(B_1^n)|^{1/n}\ls \frac{c_4}{\sqrt{n}}|T(D)|^{1/n}.\end{equation*}
It follows that
\begin{equation*}|B_1^n|^{1/n}\ls \frac{c_4}{\sqrt{n}}|D|^{1/n}.\end{equation*}
From Urysohn's inequality (see \cite{AGA-book}) we know that $w(D)\gr c_5\sqrt{n}\,|D|^{1/n}$, and a direct computation shows that
$w(B_1^n)\ls c_6\sqrt{n\log n}|B_1^n|^{1/n}$. This shows that
\begin{equation*}w(D)\gr c_7\sqrt{n/\log n}w(B_1^n).\end{equation*}
Since $D\supset B_1^n$ was arbitrary, we conclude that
\begin{equation}d_w(B_1^n)\gr c\sqrt{n/\log  n},\end{equation}
where $c>0$ is an absolute constant.

\bigbreak

\bigskip

\bigskip

%%%%%%%%%%%%%%%%%%%%%%%%%%%%%%%%%%%%%%%%%%%%%%%%%%%%%%%%%%%%%%%%%%%%%%%%%%%%%%%%%%%%%%%%%%%%%%%%%%%%%%%%%%%%%%%%%%%%%%%%%%%%%%%%%%%%%%%%%%%%%%%
\noindent {\bf Acknowledgements.} The second named author was partially supported by the US National Science Foundation grant DMS-1265155.
%%%%%%%%%%%%%%%%%%%%%%%%%%%%%%%%%%%%%%%%%%%%%%%%%%%%%%%%%%%%%%%%%%%%%%%%%%%%%%%%%%%%%%%%%%%%%%%%%%%%%%%%%%%%%%%%%%%%%%%%%%%%%%%%%%%%%%%%%%%%%%%

\bigskip

\bigskip

\footnotesize
\bibliographystyle{amsplain}

\bigskip

\bigskip

\thanks{\noindent {\bf Keywords:}  Convex bodies; Busemann-Petty problem; Shephard problem;
Sections and Projections; Volume difference inequalities; Intersection bodies; Isotropic convex body.}

\smallskip

\thanks{\noindent {\bf 2010 MSC:} Primary 52A20; Secondary 46B06, 52A23, 52A40.}

\bigskip

\bigskip

\noindent \textsc{Apostolos \ Giannopoulos}: Department of Mathematics, National and Kapodistrian University of Athens,
Panepistimiopolis 157-84, Athens, Greece.

\smallskip

\noindent \textit{E-mail:} \texttt{apgiannop@math.uoa.gr}

\bigskip

\noindent \textsc{Alexander \ Koldobsky}: Department of
Mathematics, University of Missouri, Columbia, MO 65211.

\smallskip

\noindent \textit{E-mail:} \texttt{koldobskiya@missouri.edu}

\bigskip

\end{document}